\documentclass[final,authoryear,5p,times,twocolumn]{autart}
\usepackage{amsfonts}
\usepackage{amsmath}
\usepackage{amssymb}
\usepackage{algorithm}
\usepackage{booktabs}
\usepackage{graphicx}
\usepackage{natbib}

\newtheorem{theorem}{Theorem}
\newtheorem{corollary}{Corollary}
\newtheorem{problem}{Problem}
\newtheorem{remark}{Remark}
\newtheorem{lemma}{Lemma}
\newtheorem{definition}{Definition}

\newcommand{\nx}{{m_{x}}}
\newcommand{\ny}{{m_{y}}}
\newcommand{\nt}{{m_{\theta}}}
\newcommand{\Pq}{{\mathrm{Pr}_{q}}}

\newcommand{\ev}{{V}}
\newcommand{\ic}{{k}} %iteration counter either \ell or k
\newcommand{\ms}{{{\rm\bf q}}}

\makeatletter
\renewcommand*{\@opargbegintheorem}[3]{\trivlist
      \item[\hskip \labelsep{\bfseries #1\ #2}] \textbf{(#3)}\ \itshape}
\makeatother

\begin{document}
\begin{frontmatter}

\title{A Statistical Learning Theory Approach for Uncertain Linear and Bilinear Matrix Inequalities\thanksref{footnoteinfo}}

\thanks[footnoteinfo]{This paper was not presented at any IFAC
meeting. Corresponding author Roberto Tempo. Tel. +39 011 090-5408.
Fax +39 011 090-5429.}

\author[DSI,NUS]{Mohammadreza Chamanbaz}\ead{Mohammad\_C@dsi.A-Star.edu.sg and Chamanbaz@nus.edu.sg},
\author[IEIIT]{Fabrizio Dabbene}\ead{fabrizio.dabbene@polito.it},
\author[IEIIT]{Roberto Tempo}\ead{roberto.tempo@polito.it},
\author[DSI]{Venkatakrishnan Venkataramanan}\ead{Venka\_V@dsi.A-Star.edu.sg},
\author[NUS]{Qing-Guo Wang}\ead{elewqg@nus.edu.sg}

\address[DSI]{Data Storage Institute, Singapore.}
\address[IEIIT]{CNR-IEIIT, Politecnico di Torino, Torino 10129, Italy.}
\address[NUS]{Department of Electrical and Computer Engineering, National University of Singapore.}

\begin{abstract}
In this paper, we consider the problem of minimizing  a linear functional subject to uncertain linear and bilinear matrix inequalities, which  depend in a possibly nonlinear way on a vector of uncertain parameters.
Motivated by recent results in statistical learning theory, we show that probabilistic guaranteed solutions can be obtained by means of randomized algorithms. In particular, we show that the Vapnik-Chervonenkis dimension (VC-dimension) of the two problems is finite, and we compute  upper bounds on it. In turn, these bounds allow us to derive explicitly the sample complexity of these problems. Using these bounds, in the second part of the paper, we derive a sequential scheme, based on a sequence of optimization and validation steps. The algorithm is on the same lines of recent schemes proposed for similar problems, but improves both in terms of complexity and generality. The effectiveness of this approach is shown using a linear model of a robot manipulator subject to uncertain parameters.
\end{abstract}


\end{frontmatter}

\section{Introduction}\label{sec:intro}

Statistical learning theory is a very effective tool in dealing with various applications, which include neural networks and control systems, see for instance \cite{vidyasagar_learning_2002}.
The main objective of this theory is to extend convergence properties of the empirical mean, which can be
computed with a Monte Carlo simulation, from finite families to infinite families of functions. For finite families, these properties can be easily established by means of a repeated application of the so-called Hoeffding inequality, and are related to the well-known law of large numbers, see for instance \cite{tempo_randomized_2012}. On the other hand, for infinite families  deeper technical tools have been developed in the seminal work of \cite{VapChe:71}. In this case, the main issue is to establish \textit{uniform} convergence of empirical means. In particular, this requires to determine whether or not a combinatorial parameter called the Vapnik-Chervonenkis dimension (VC-dimension) is finite, see \cite{vapnik_statistical_1998}.

Subsequent contributions on statistical learning theory by \cite{vidyasagar_randomized_2001}
followed two main research directions: First, to demonstrate that this theory is indeed an effective tool for control of systems affected by uncertainty. Second, to ``invert" the bounds provided by Vapnik and Chervonenkis, introducing the concept of \textit{sample complexity}. Roughly speaking, when dealing with control of uncertain systems, the sample complexity provides the number of random samples of the uncertainty that should be drawn to derive a stabilizing controller (or a controller which attains a given $\mathcal{H}_\infty$-norm bound on the closed-loop sensitivity function), with sufficiently high probabilistic accuracy and confidence. Since the sample complexity is a function of the accuracy, confidence and the VC-dimension, specific bounds
on this combinatorial parameter should be derived.
In turn, this involves a problem reformulation in terms of Boolean functions, and the evaluation of the number of required polynomial inequalities, their order and the number of design variables.

For various stabilization problems, which include stability of interval matrices and simultaneous stabilization with static output feedback, bounds on the VC-dimension have been derived  in \cite{VidBlo:01}. In this paper, we continue this specific line of research, and we compute the VC-dimension for control problems formulated in terms of \textit{uncertain} linear matrix inequalities (LMIs) and bilinear matrix inequalities (BMIs). It is well-known that many robust and optimal control problems can be indeed formulated in these forms, see for instance \cite{boyd_linear_1994,GoSaPa:95,kanev_robust_2004,vanantwerp_tutorial_2000}. 
When the uncertain LMIs depend on the uncertainty in a linear or multilinear way, extreme point results can be derived and applied, see e.g.~\cite{ATRC:08,CalDab:08b}.
For more general linear fractional dependence on norm-bounded \emph{unstructured}  uncertainty, efficient methods based on the so-called robustness lemma have been successfully developed, see e.g.~\cite{ElGLeb:98}.
However, the problem is significantly harder to solve when it involves structured uncertainty and computable solutions can be in general obtained only at the expense of introducing conservatism.

The main contribution of the present paper is to compute upper bounds on the VC-dimension of uncertain LMIs and BMIs, and
to establish the related sample complexity. 
We remark that the sample complexity is independent from the number of uncertain parameters entering into the LMIs and BMIs, 
and on their functional relationship. Hence, the related randomized algorithms run
in polynomial-time. However, for relatively small values of the
probabilistic accuracy and confidence, the sample complexity turns out to be very large, as usual in the context of statistical learning theory.
For this reason, randomized algorithms based on a direct application of these bounds may be of limited use in practice.
To alleviate this difficulty, in the second part of the paper we propose a sequential algorithm specifically tailored to the problem at hand. This algorithm has some similarities with  sequential algorithms developed for other problems in the area of randomized algorithms for control of uncertain systems, see
Remark \ref{rem:related} for a discussion of these results.

Finally, the effectiveness of this approach is shown by a numerical example related to the static output feedback stabilization of an uncertain robot manipulator joint taken from
 \citep{kanev_controller_2000}.
In particular, the objective is to design a static output feedback controller which minimizes the worst-case $\mathcal{H}_\infty$ norm. The numerical performance of the proposed sequential algorithm is evaluated and compared with the theoretical sample-complexity previously derived.

\section{Problem Formulation}\label{sec:problem formulation}
%Many robust and optimal control problems can be formulated as linear or bilinear matrix inequality (LMI or BMI). In the case where problem data involves  uncertain parameters the LMI and BMI problems are in the form of semi-infinite optimization programs, due to the infinite number of constraints involved. 
We now formally state the uncertain LMI and BMI problems discussed in the Introduction.

\begin{problem}[Uncertain strict LMI optimization]\label{prob:LMI}
Find the optimal value of $x$, if it exists, which solves the optimization problem
\begin{align} \label{eq:LMI}
& \min_{x}\, c_x^Tx\quad \text{\rm subject to}\\
&\quad F_\textrm{LMI}(x,q)\doteq F_0(q)+\sum_{i=1}^{\nx} x_iF_i(q)\succ0,\,\forall q\in\mathbb{Q}
\nonumber
\end{align}
where $x\in\mathbb{R}^\nx$ is the vector of optimization variables, $q\in\mathbb{Q}\subset \mathbb{R}^\ell$ is the
vector of uncertain parameters bounded in the set $\mathbb{Q}$ and $F_i=F_i^T\in\mathbb{R}^{n\times n},\,i=0,\ldots,\nx$. The inequality $F_\textrm{LMI}(x,q)\succ0$ means that $F_\textrm{LMI}(x,q)$ is positive definite.% , i.e. $u^TF_\textrm{LMI}(x,q)u>0\text{ for all }u\in\mathbb{R}^n$.
\end{problem}

\begin{problem}[Uncertain strict BMI optimization]\label{prob:BMI}
Find the optimal values of $x$ and $y$, if they exist, which solve the optimization problem
\begin{align}
& \min_{x,y}\, c_x^Tx+c_y^Ty\quad \text{\rm subject to}  \label{eq:BMI} \\
&\quad F_\textrm{BMI}(x,y,q)\doteq F_0(q)+\sum_{i=1}^\nx x_iF_i(q)+\sum_{j=1}^\ny y_jG_j(q)\nonumber\\
& \qquad\qquad\qquad+\sum_{i=1}^\nx \sum_{j=1}^\ny x_iy_jH_{ij}(q)\succ0, \,\forall q\in\mathbb{Q} \nonumber
\end{align}
where $x\in\mathbb{R}^\nx$ and $y\in\mathbb{R}^\ny$ are the vectors of optimization variables, $q\in\mathbb{Q}\subset\mathbb{R}^\ell$ is the vector of uncertain parameters,  $F_0=F_0^T\in\mathbb{R}^{n\times n}$, and $F_i=F_i^T\in\mathbb{R}^{n\times n}$, $G_j=G_j^T\in\mathbb{R}^{n\times n}$, $H_{ij}=H_{ij}^T\in\mathbb{R}^{n\times n}$, $i=1,\ldots,\nx$, $j=1,\ldots,\ny$.
%The inequality $F_\textrm{BMI}(x,q)\succ0$ has the same interpretation as in Problem 1.
\end{problem}

In order to allow a unified treatment of Problems 1 and 2, we now formally 
define the design parameters for LMIs and BMIs.

\begin{definition}[Design parameters for LMIs/BMIs]\label{def:LMI-BMI}
For Problem 1, we define
\[
\theta \doteq x, \quad \nt=\nx,\quad c_{\theta}=c_{x}; \label{def-LMI}
\]
and, for Problem 2, we define
\[
\theta \doteq \displaystyle{x\brack y},\quad \nt={\nx+\ny}, \quad\displaystyle c_{\theta}={c_x \brack c_y}. \label{def-BMI}
\]
\end{definition}

%In the present paper, in the spirit of \cite{tempo_randomized_2012}, we study a probabilistic framework for solving Problems \ref{prob:LMI} and \ref{prob:BMI}. Formally, 
Next, we assume that $q$ is a random variable and  a probability measure $\Pq$ over the Borel $\sigma$-algebra of $\mathbb{Q}\subset\mathbb{R}^\ell$ is given. 
Then, the constraints in (\ref{eq:LMI}) and (\ref{eq:BMI}) become chance-constraints, see e.g.\ \cite{Uryasev:00}, which may be violated for some $q\in\mathbb{Q}$. This concept is formally expressed using the notion of ``probability of violation".

\begin{definition}[Probability of violation]\label{def:prob of viol}
%Consider the notation introduced in (\ref{def-LMI})-(\ref{def-BMI}). 
The probability of violation of $\theta$ for the binary-valued function $g:\mathbb{R}^\nt\times\mathbb{Q}\to\{0,1\}$ is defined as
\begin{equation}\label{eq:prob of viol}
V_g(\theta)\doteq \Pq\left\{q\in\mathbb{Q}\,\,:\,g(\theta,q)=1\right\}
\end{equation}
where, for  Problem \ref{prob:LMI}, 
\begin{equation}\label{eq:g for lmi}
g(\theta,q)\doteq\begin{cases}
0\quad \text{if}\,F_\textrm{LMI}(\theta,q)\succ0\\
1\quad \text{otherwise}
\end{cases}
\end{equation}
and, for  Problem \ref{prob:BMI}, 
\begin{equation}\label{eq:g for bmi}
g(\theta,q)\doteq\begin{cases}
0\quad \text{if}\,F_\textrm{BMI}(\theta,q)\succ0\\
1\quad \text{otherwise}
\end{cases}.
\end{equation}
\end{definition}

We remark that the probability  of violation is in general very hard to evaluate, due to the difficulty of computing  the multiple integrals associated with (\ref{eq:prob of viol}). Nevertheless, we can ``estimate" this probability using randomization. To this end, we extract $N$ independent  identically distributed (i.i.d) samples from the set $\mathbb{Q}$
\[
\ms=\{q^{(1)},\ldots,q^{(N)}\}\in\mathbb{Q}^{N},
\]
according to the measure $\Pq$, where $\mathbb{Q}^{N}\doteq \mathbb{Q}\times\mathbb{Q}\times\cdots\times\mathbb{Q}$ ($N$ times). Next, a Monte Carlo approach is employed to obtain the so called ``empirical violation";
see e.g.\ \cite{vidyasagar_learning_2002}.
\begin{definition}[Empirical violation]\label{def:emp viol}
For given $\theta\in\mathbb{R}^\nt$ the empirical violation of $g(\theta,q)$ with respect to the multisample $\ms=\{q^{(1)},\ldots,q^{(N)}\}$ is defined as
\begin{equation}\label{eq:emp mean}
\widehat{V}_g(\theta,\ms)\doteq\frac{1}{N}\sum _{i=1}^N g(\theta,q^{(i)}).
\end{equation}
\end{definition}

\subsection{Randomized Strategy to Optimization Problems}
There are several randomized methodologies in the literature which are based on randomization in the uncertainty space, design parameter space or both. For example, in  \cite{vidyasagar_randomized_2001} randomization in both uncertainty and design parameter spaces is employed for minimizing the empirical mean. Similarly, a min-max approach and a bootstrap learning method  are presented in \cite{calafiore_probabilistic_2006}
and \cite{koltchinskii_improved_2000}, respectively, but
 these papers deal with finite families. In \cite{alamo_randomized_2009} the authors proposed a randomized algorithm for infinite families which is applicable to convex and non-convex problems. Finally, a non-sequential randomized methodology for uncertain convex problems is introduced in \cite{calafiore_uncertain_2004,calafiore_scenario_2006,campi_exact_2008}.
In  Algorithm \ref{alg:randomized algorithm} we present a non-sequential  randomized strategy for solving Problems \ref{prob:LMI} and \ref{prob:BMI}.

\begin{algorithm}[h!]
\caption{\textsc{A Randomized Strategy for Uncertain LMIs/BMIs}}
\label{alg:randomized algorithm}
\begin{itemize}
  \item Given the underlying probability density function (pdf) over the uncertainty set $\mathbb{Q}$ and the level parameter $\rho\in[0,1)$, extract $N$ independent identically distributed  samples from $\mathbb{Q}$ based on the underlying pdf
      \[
      \ms=\{q^{(1)},\ldots,q^{(N)}\}.
      \]
  \item Find the optimal value, if it exists, of the following optimization problem
  \begin{align}\label{eq:samples optimization} \nonumber
  & \underset{\theta}{\text{minimize}}\quad c_{\theta}^T\theta\\
  & \text{subject to}\quad \widehat{V}_g(\theta,\ms)\leq\rho
  \end{align}
\end{itemize}
\end{algorithm}
We remark that introducing the level parameter $\rho>0$ enables us to handle probabilistic (soft) constraints,
in the same spirit of \cite{alamo_randomized_2009}.
The main objective of the present paper is to derive the explicit sample complexity bound on $N$ based on statistical learning theory results. 
Finally, we remark that, in the case of LMI constraints, problem (\ref{eq:samples optimization}) is a semidefinite optimization problem (SDP) that  can be solved efficiently, see  
\cite{VanBoy:96} and \cite{Todd:01} for a discussion on the numerical aspects of solving SDP problems. In the case of BMI constraints, 
efficient solvers such as PENBMI \citep{kocvara_pennon:_2003} are available, but global solutions to the optimization problem in general cannot be  obtained.

\section{Vapnik-Chervonenkis Theory}\label{sec:VC theory}
In this section, we give a very brief overview of the Vapnik-Chervonenkis theory. The material presented is classical, but a summary is instrumental to our next developments. In particular, we review some bounding inequalities which are  used in the subsequent sections to derive the explicit sample bounds for solving Problems \ref{prob:LMI} and \ref{prob:BMI}.

\begin{definition}[Probability of two-sided failure]\label{def: prob two fail}
Given $N,\,\varepsilon\in(0,1)$ and $g:\mathbb{R}^\nt\times\mathbb{Q}\to\{0,1\}$, the probability of two-sided failure denoted by $q_g(N,\varepsilon)$ is defined as
\[
q_g(N,\varepsilon)\doteq\Pq\left\{\ms\in\mathbb{Q}^N: \underset{\theta\in\mathbb{R}^\nt}{\sup}|V_g(\theta)-\widehat{V}_g(\theta,\ms)|>\varepsilon\right\}.
\]
\end{definition}

The probability of two-sided failure determines how close the empirical violation is to the true probability of violation. In other words, if we extract a multisample $\ms$ with cardinality $N$ from the uncertainty set $\mathbb{Q}$, we guarantee that the empirical violation (\ref{eq:emp mean}) is within $\varepsilon$ of the true probability of violation (\ref{eq:prob of viol}) for all $q\in\mathbb{Q}$ except for a subset having probability measure at most $q_g(N,\varepsilon)$. The parameter $\varepsilon\in(0,1)$ is called accuracy.

Let $\mathcal{G}$ denote the family of  functions $\{g(\theta,q):\theta\in\mathbb{R}^\nt\}$ where $g:\mathbb{R}^\nt\times\mathbb{Q}\to\{0,1\}$ is defined in (\ref{eq:g for lmi}) or in (\ref{eq:g for bmi}). The family $\mathcal{G}$ is said to satisfy the property of \emph{uniform convergence of empirical mean} (UCEM) if $q_g(N,\varepsilon)\to 0\text{ as }N\to\infty$ for any $\varepsilon\in(0,1)$. We remark that if $\mathcal{G}$ includes  finite family of functions, it indeed has the UCEM property. However, infinite families do not necessarily enjoy the UCEM property, see \cite{vidyasagar_learning_2002} for several examples of this type. Problems \ref{prob:LMI} and \ref{prob:BMI} belong to the class of infinite family of  functions.

We define the family $\mathcal{S}_g$ containing all possible sets
$
S_g\doteq\left\{q\in\mathbb{Q}: g(\theta,q)=1\right\}
$,
for $g$ varying in $\mathcal{G}$. Now consider a multisample $\ms=\{q^{(1)},\ldots,q^{(N)}\}$ of cardinality $N$. For the family of  functions $\mathcal{G}$, let
\[
\mathbb{N}_\mathcal{G}(\ms)\doteq\text{Card}\left(\ms\cap S_g,S_g\in\mathcal{S}_g\right).
\]
In words, we say that $\mathcal{S}_g$ ``shatters'' $\ms$ when $\mathbb{N}_\mathcal{G}$ is equal to $2^N$. The notion of ``shatter coefficient", also known as ``growth function", is now defined formally.
\begin{definition}[Shatter Coefficient]
The shatter coefficient of the family $\mathcal{G}$, denoted by $\mathbb{S}_\mathcal{G}(N)$, is defined as
\[
\mathbb{S}_\mathcal{G}(N)\doteq\underset{\ms\in\mathbb{Q}^N}{\max}\mathbb{N}_\mathcal{G}(\ms).
\]
\end{definition}
A bound on the shatter coefficient can be obtained by Sauer lemma \citep{sauer_density_1972}, which in turn requires the computation of the VC-dimension, defined next.
\begin{definition}[VC-dimension]
The {\rm VC}-dimension of the family of  functions $\mathcal{G}$ is defined as the largest integer $d$ for which $\mathbb{S}_\mathcal{G}(N)=2^d$.
\end{definition}
The following result establishes a bound on the probability of two-sided failure in terms of VC-dimension.
\begin{theorem}[Vapnik and Chervonenkis]\label{theorem:vc result}
Let $d$ denote the {\rm VC}-dimension of the family of  functions $\mathcal{G}$. Then, for any $\varepsilon\in(0,1)$
\begin{equation}\label{eq:vc result}
q_g(N,\varepsilon)\leq4e^{2\varepsilon}\left(\frac{2eN}{d}\right)^d e^{-N\varepsilon^2}
\end{equation}
where $e$ is the Euler number.
\end{theorem}
This result has been proven in \cite{VapChe:71} and it is stated in \cite[Theorem 4.4]{vapnik_statistical_1998}.
%see also  \cite[Corollary 1]{alamo_randomized_2009}.

\section{Main Results}\label{sec: main results}
In  view of Theorem \ref{theorem:vc result}, we conclude that families with finite VC-dimension $d<\infty$ enjoy the UCEM property. Hence, it is important i) to show that the collection $\mathcal{G}$ of  functions has finite VC-dimension and, ii) to derive upper bounds on the VC-dimension.

\subsection{Computation of Vapnik-Chervonenkis Dimension}\label{sec:computation of VC dim}
In the next theorem, which is one of the main contributions of this paper, we derive an upper bound on the VC-dimension of the uncertain LMI and BMI in Problems \ref{prob:LMI} and \ref{prob:BMI}.

\begin{theorem}[VC bounds for strict LMIs/BMIs]\label{theorem:VC dimension of LMI and BMI}
Consider the notation introduced in Definition \ref{def:LMI-BMI}. 
Then, the {\rm VC}-dimension of uncertain LMIs and BMIs (Problems \ref{eq:LMI} and \ref{eq:BMI}) is upper bounded by 
\begin{equation}
\label{th1-MI}
d \le 2\nt\lg (4en^2)
\end{equation}
where $\lg(.)$ denotes the logarithm to the base~$2$.
\end{theorem}
\noindent
\textbf{Proof}
See Appendix \ref{app:proof of vc dim for lmi and bmi}.

It is interesting to observe that the VC-dimension of uncertain LMIs and BMIs is linear in the number of design variables $\nt$. 
%We also remark that the VC-dimension for the general case of nonlinear matrix inequality (NMI) cannot be computed. This is due to the fact that there is no optimization %variable of degree larger than one in the BMI, and this  is clearly not the case for NMI. 
In the next subsection, we derive explicit sample bounds to be used in Algorithm \ref{alg:randomized algorithm} for solving Problems \ref{prob:LMI} and \ref{prob:BMI}.

\subsection{Sample Complexity Bounds}\label{sec:sample bounds}
In  this section, we study a number of sample bounds guaranteeing that the probability of failures is bounded by a confidence parameter $\delta\in(0,1)$.
We remark that there are several results in the literature to derive sample complexity bounds. To the best of our knowledge, the least conservative is stated in
 Corollary 3 in \cite{alamo_randomized_2009}. For given $\varepsilon,\delta\in(0,1)$, the probability of two-sided failure (see Definition \ref{def: prob two fail}) is bounded by $\delta$ provided that at least
\begin{equation}\label{eq:two sided sample bound}
N\geq\frac{1.2}{\varepsilon^2}\left(\ln\frac{4e^{2\varepsilon}}{\delta}+d\ln\frac{12}{\varepsilon^2}\right)
\end{equation}
samples are drawn, where $d<\infty$ denotes the VC-dimension of the family of  functions $\mathcal{G}$,  and $\ln$ is the natural logarithm.
This result is exploited in the next corollary, that provides the explicit sample complexity bound for the probability of two-sided failure.

\begin{corollary}\label{corollary:two sided sample bound}
Consider the notation introduced in Definition \ref{def:LMI-BMI}, and suppose that $\varepsilon, \delta\in(0,1)$ are given. Then, the probability of two-sided failure is bounded by $\delta$ if at least
%\begin{equation}\label{eq:sample bound LMI twosided}
%N_{{\rm LMI}}\geq\frac{1.2}{\varepsilon^2}\left(\ln\frac{4e^{2\varepsilon}}{\delta}+2\nx\lg(4en^2)\ln\frac{12}{\varepsilon^2}\right)
%\end{equation}
%and
\begin{equation}\label{eq:sample bound BMI twosided}
N\geq\frac{1.2}{\varepsilon^2}\left(\ln\frac{4e^{2\varepsilon}}{\delta}+2\nt\lg(4en^2)\ln\frac{12}{\varepsilon^2}\right)
\end{equation}
samples are drawn for the Problems \ref{prob:LMI} and \ref{prob:BMI}.
\end{corollary}
\noindent
\textbf{Proof}
The statement of Corollary  \ref{corollary:two sided sample bound} follows immediately by combining (\ref{eq:two sided sample bound}) and the results of Theorem \ref{theorem:VC dimension of LMI and BMI}.

A weaker notion than the probability of two-sided failure is the ``probability of one-sided constrained failure" introduced in the following definition.
\begin{definition}[Probability of one-sided constrained failure]
Given $N,\,\varepsilon\in(0,1),\,\rho\in[0,1)$ and $g:\mathbb{R}^\nt\times\mathbb{Q}\to\{0,1\}$, the probability of one-sided constrained failure, denoted by $p_g(N,\varepsilon,\rho)$, is defined as
\begin{eqnarray*}
\lefteqn{p_g  (N,\varepsilon,\rho)\doteq  \Pq\Bigg\{\ms\in\mathbb{Q}^N:\text{ there exists }\theta\in\mathbb{R}^\nt}\\
 && \quad\text{ such that }\widehat{V}_g(\theta,\ms)\leq\rho \text{ and }V_g(\theta)-\widehat{V}_g(\theta,\ms)>\varepsilon\Bigg\}.
\end{eqnarray*}
\end{definition}
Following the same lines of Corollary \ref{corollary:two sided sample bound}, sample complexity bounds for the probability of one-sided constrained failure are derived.

\begin{corollary}\label{corollary:one sided sample bound}
Consider the notation introduced in Definition \ref{def:LMI-BMI}, and suppose that $\varepsilon\in(0,1)$, $\delta\in(0,1)$ and $\rho\in[0,1)$ are given. Then,  the probability of one-sided constrained failure is bounded by $\delta$ if at least
\begin{equation}\label{eq:sample bound BMI onesided}
N\geq\frac{5(\rho+\varepsilon)}{\varepsilon^2}\left(\ln\frac{4}{\delta}+2\nt\lg(4en^2)\ln\frac{40(\rho+\varepsilon)}{\varepsilon^2}\right)
\end{equation}
samples are drawn for the Problems \ref{prob:LMI} and \ref{prob:BMI}.
\end{corollary}

\noindent
\textbf{Proof}
This result is an immediate consequence of Theorem 7 in \cite{alamo_randomized_2009}, which
states that, for given $\varepsilon,\delta\in(0,1)$ and $\rho\in[0,1)$, the probability of one-sided constrained failure is bounded by $\delta$ provided that at least
\begin{equation}\label{eq:one sided sample bound}
N\geq\frac{5(\rho+\varepsilon)}{\varepsilon^2}\left(\ln\frac{4}{\delta}+d\ln\frac{40(\rho+\varepsilon)}{\varepsilon^2}\right)
\end{equation}
samples are drawn, where $d<\infty$ denotes the VC-dimension of the family of  functions $\mathcal{G}$.
The statement in Corollary  \ref{corollary:one sided sample bound} is derived by substituting the results of Theorem \ref{theorem:VC dimension of LMI and BMI} into (\ref{eq:one sided sample bound}).

\begin{figure}[h]
\centerline{\includegraphics[width=9cm]{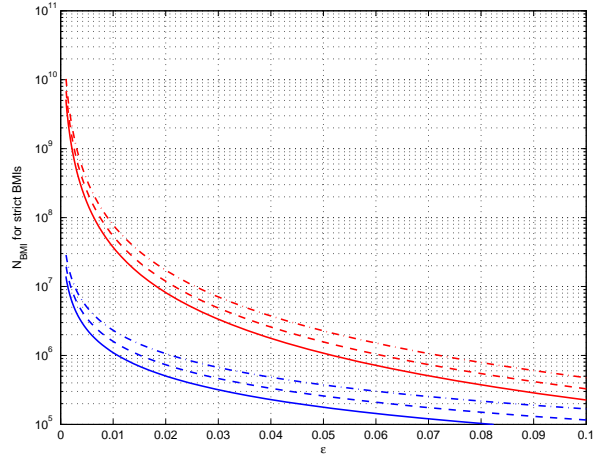}}
\caption{Sample complexity bounds for strict BMIs,
for $\delta=1\times10^{-8}$, $m_{x}+m_{y}=13$, and for different BMI dimensions:
$n=10$ (continuous line) $n=50$ (dashed line) and $n=100$ (dash-dotted line). The red (upper) plots show the two-sided bound (\ref{eq:sample bound BMI twosided}),
while the blue (lower) plots show the one-sided constrained failure bound (\ref{eq:sample bound BMI onesided}) for $\rho=0$.}
\label{fig:sample bounds strict}
\end{figure}

Note that the sample complexity of Corollary \ref{corollary:one sided sample bound} improves upon that of
 Corollary~\ref{corollary:two sided sample bound}, as shown in Figure~\ref{fig:sample bounds strict}. In particular, it is clear that the bound (\ref{eq:sample bound BMI twosided}) grows as $\mathcal{O}(\frac{1}{\varepsilon^2}\ln\frac{1}{\varepsilon^2})$, which implies that if the accuracy level $\varepsilon$ is chosen to be very small, the sample bounds  can be very large, while
 (\ref{eq:sample bound BMI onesided}) grows as $\mathcal{O}(\frac{1}{\varepsilon}\ln\frac{1}{\varepsilon})$.

\section{Semidefinite Constraints}\label{sec:semidefinite constraints}
In this section, we compute upper bounds on the VC-dimension of the semidefinite versions of Problems \ref{prob:LMI} and~\ref{prob:BMI} where strict inequalities $(\succ0)$ are replaced with nonstrict inequalities $(\succcurlyeq0)$%
\footnote{Throughout the paper, nonstrict (semidefinite) versions of Problems \ref{prob:LMI} and \ref{prob:BMI} are called ``uncertain semidefinite LMI  problem" and ``uncertain semidefinite BMI  problem" respectively.}.
 Semidefinite constraints appear in some control problems such as dissipativity; furthermore, some modeling languages such as YALMIP \citep{YALMIP} treat strict inequalities using  nonstrict ones by adding a slight perturbation. Hence, it is important to derive sample complexities for uncertain semidefinite LMI and BMI problems.

In the following theorem, we establish upper bounds on the VC-dimension of uncertain semidefinite LMI and BMI  problems. The proof 
of this result is reported in Appendix \ref{app:proof of nonstrict}.

\begin{theorem}[VC bounds for nonstrict LMIs/BMIs]\label{theorem:VC dim semidef}
Consider the notation introduced in Definition \ref{def:LMI-BMI}. Then, the {\rm VC}-dimension of uncertain semidefinite LMI and  BMI problems  is upper bounded by 
\[
d\le2\nt \lg (4en2^n).
\]
\end{theorem}

\begin{remark}[Strict and nonstrict LMIs/BMIs]
Comparing the bounds of Theorems \ref{theorem:VC dimension of LMI and BMI}
and \ref{theorem:VC dim semidef}, it can be seen that the bounds on the VC-dimension
of strict and nonstrict  LMIs/BMIs differ only in the terms $n^{2}$ and $n2^n$ appearing in the arguments of the logarithm. That is, the quadratic dependence on $n$ of strict LMIs/BMIs becomes exponential
for nonstrict ones. Note however that this effect is largely mitigated by the logarithm.
This difference is not surprising, and it follows from the fact that checking positive semi-definiteness requires non-negativity of all principle minors, as discussed in Appendix \ref{app:proof of nonstrict}. 
To show this fact, consider the matrix
{\small
\vspace{-6mm}
\[
\left[\begin{array}{ccc}
1 & 1 & 1 \\
1 & 1 & 1\\
1 & 1 & 0 \\
\end{array}\right]
\vspace{-6mm}
\]}

\noindent
introduced in \cite{bernstein2009matrix}. This matrix 
has leading principal minors equal to $1$, $0$ and $0$, which are nonnegative, but it is  not positive semidefinite, because its eigenvalues are $2.732$, $0$, and $-0.732$. Note that the same issue arises in \cite[Theorem 4]{VidBlo:01}, regarding positive definiteness and semi-definiteness of interval matrices.
\end{remark}

\begin{figure}[h]
\centerline{\includegraphics[width=9cm]{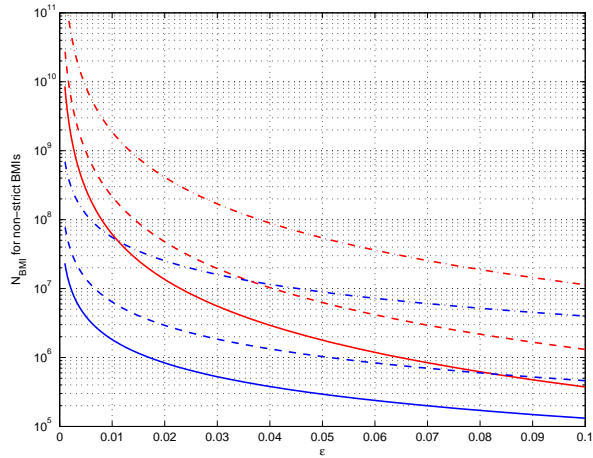}}
\caption{Sample complexity bounds for nonstrict BMIs,
for $\delta=1\times10^{-8}$, $m_{x}+m_{y}=13$, and for different BMI dimensions:
$n=10$ (continuous line) $n=50$ (dashed line) and $n=100$ (dash-dotted line). The red (upper) plots show the two-sided bound,
while the blue (lower) plots show the one-sided constrained failure bound  for $\rho=0$.}
\label{fig:sample bounds-NS}
\end{figure}

\begin{remark}[Explicit sample complexity for nonstrict LMIs/BMIs]\label{rem:nonstrict}
Using the results of Theorem \ref{theorem:VC dim semidef}, we can establish bounds on sample complexity which guarantee the probability of two-sided failure and the probability of one-sided constrained failure of uncertain semidefinite LMI and BMI problems to be bounded by the confidence parameter $\delta$. It should be noted that for semidefinite problems of this section, Definition \ref{def:prob of viol} is revised accordingly such that  strict inequalities  in (\ref{eq:g for lmi}) and (\ref{eq:g for bmi}) are replaced  with nonstrict ones. 
This also affects empirical violation, probability of two-sided failure and probability of one-sided constrained failure. 

Then, the results of Corollaries  \ref{corollary:two sided sample bound} and \ref{corollary:one sided sample bound} for the uncertain  semidefinite LMI and BMI problems immediately hold provided that the VC-dimension bound  $2\nt \lg(4en^2)$ is replaced by $2\nt \lg(4en2^n)$.
The sample complexity bounds for semidefinite BMIs are illustrated in Figure~\ref{fig:sample bounds-NS}.
\end{remark}

It should be also noted that the sample complexity bounds derived in this paper for the uncertain strict and semidefinite LMI and BMI problems can be quite large. This is a usual situation in the context of statistical learning, that may lead to computationally expensive optimization problems if Algorithm \ref{alg:randomized algorithm} is applied in one-shot. This motivates the developments of the next section, where a \textit{sequential} randomized algorithm for bounding the probability of one-sided constrained failure is presented. The sequential algorithm can alleviate the computational burden of directly solving (\ref{eq:samples optimization}).

\section{Sequential Randomized Algorithm}\label{sec:sequential algorithm}

Sequential methods in probabilistic design usually follow an iterative scheme which includes optimization steps to update the design parameters, followed by randomization steps to check the feasibility of the candidate solution \citep{tempo_randomized_2012}. The first step is deterministic, while the second one is probabilistic. Examples of such scheme are probabilistic design methods based on gradient \citep{PolTem:01,CalPol:01}, ellipsoid \citep{KaDeVe:03,Oishi:07} and cutting plane \citep{calafiore_probabilistic_2007} update rules, see \cite{tempo_randomized_2012} for more details.

Recently \cite{alamo_randomized_2012} introduced a  general framework for nonconvex problems, defining the class of sequential probabilistic validation (SPV) algorithms.
In this section we propose a sequential randomized algorithm specifically tailored for the problem at hand, which mitigates the conservatism of the bound (\ref{eq:sample bound BMI onesided}) or its corresponding sample bound for the uncertain semidefinite LMI and BMI problems. This is accomplished by generating a sequence of ``design" sample sets $\{q_d^{(1)},\ldots,q_d^{(N_\ic)}\}$ with increasing cardinality $N_\ic$ which are used in (\ref{eq:samples optimization}) for solving the optimization problem. In parallel, ``validation" sample sets  $\{q_v^{(1)},\ldots,q_v^{(M_\ic)}\}$ of cardinality $M_\ic$ are also generated by the algorithm in order to check whether the given candidate solution, obtained from solving (\ref{eq:samples optimization}), satisfies the desired  probability of violation. The proposed scheme is reported in Algorithm \ref{alg:sequential algorithm}.

%For simplicity of notation, we denote the  sample bounds derived in  Corollary  \ref{corollary:one sided sample bound} and the corresponding sample complexity for the uncertain semidefinite LMI and BMI problems by $N_{\rm MI}$.

%%%%%%%%
\begin{algorithm}[h!]
\caption{\textsc{A Sequential Randomized Algorithm}}
\label{alg:sequential algorithm}
\begin{enumerate}
  \item \textsc{Initialization} \newline
  Set the iteration counter to zero $(\ic=0)$. Choose the desired accuracy $\varepsilon\in(0,1)$, confidence $\delta\in(0,1)$ and level $\rho\in[0,1)$ parameters and the desired number of iterations $\ic_t>1$.

  \item \textsc{Update}\label{item:update}\newline
  Set $\ic=\ic+1$ and $N_\ic\ge N\frac{\ic}{\ic_t}$ where $N$ satisfies (\ref{eq:sample bound BMI onesided}).
  \item \textsc{Design}
  \begin{itemize}
    \item Draw $N_\ic$ i.i.d samples $\ms_d=\{q_d^{(1)}\ldots q_d^{(N_\ic)}\}$ from the uncertainty set $\mathbb{Q}$ based on the underlying distribution.
    \item Solve the following optimization problem
    \begin{align}\label{eq:optimization in seq alg} \nonumber
    & \underset{\theta}{\text{minimize}}\quad c_{\theta}^T\theta\\
    & \text{subject to}\quad \widehat{V}_g(\theta,\ms_d)\leq\rho.
    \end{align}
    \item \textbf{If} the optimization problem (\ref{eq:optimization in seq alg}) is not feasible, the original problem is not feasible as well.
    \item \textbf{Else if}, the last iteration is reached $(\ic=\ic_t)$, $\widehat{\theta}_{N_\ic}$ is a probabilistic robust solution and \textbf{Exit}.
   \item \textbf{Else}, continue to the next step.
  \end{itemize}
  \item \textsc{Validation}
  \begin{itemize}
    \item Draw
    \begin{equation}\label{eq:sample bound Mk}
    M_\ic\ge \frac{\alpha\ln k+\ln \left(\mathcal{S}_{k_t}(\alpha)\right)+\ln\frac{1}{\delta}}{\ln\left(\frac{1}{(\rho+\varepsilon)a^{\rho-1}+a^\rho(1-(\rho+\varepsilon))}\right)}
    \end{equation}
    i.i.d. samples $\ms_v=\{q_v^{(1)}\ldots q_v^{(M_\ic)}\}$ from the uncertainty set $\mathbb{Q}$ based on the underlying distribution.     \item \textbf{If}
        \[
        \widehat{V}_g(\widehat{\theta}_{N_\ic},\ms_v)\leq\rho
        \]
        then, $\widehat{\theta}_{N_\ic}$ is a probabilistic  solution and \textbf{Exit}.
   \item \textbf{Else}, goto step (\ref{item:update}).
  \end{itemize}
\end{enumerate}
\end{algorithm}
%%%%%%%%

Note that step (3) of this Algorithm is for the case of strict LMIs/BMIs. In the nonstrict case the bound  (\ref{eq:sample bound BMI onesided}) needs to be replaced by the bound discussed in Remark~\ref{rem:nonstrict}.
Note also that, the validation bound (\ref{eq:sample bound Mk}) in step (4), the parameters $a\geq1$ and $\alpha>0$  are real and $\mathcal{S}_{k_t}(\alpha)$ is a finite hyperharmonic series also known as $p$-series, i.e.
    \[
    \mathcal{S}_{\ic_t}(\alpha)=\sum_{\ic=1}^{\ic_t}\frac{1}{k^\alpha}.
    \]

The theoretical properties of Algorithm \ref{alg:sequential algorithm} are summarized in the next theorem, see Theorem 5 in \cite{chamanbaz_sample_2013} for proof.
\begin{lemma}\label{theo:property of seq algorithm}
Suppose $\epsilon, \delta\in(0,\,1)$ are given. Then, if at iteration $\ic$  Algorithm \ref{alg:sequential algorithm} exits with a probabilistic  solution $\widehat{\theta}_{N_\ic}$, then it holds that $\ev_g(\widehat{\theta}_{N_\ic})\leq\rho+\varepsilon$ with probability no smaller than $1-\delta$, i.e.
\[
\Pq\left\{V_g(\widehat{\theta}_{N_k})\leq\rho+\varepsilon\right\}\geq1-\delta.
\]
\end{lemma}

%%%%%%%%%

\begin{remark}[Comments on Algorithm  \ref{alg:sequential algorithm} and related results]
\label{rem:related}
Algorithm \ref{alg:sequential algorithm} follows the general scheme of other sequential algorithms previously developed  in the area of randomized algorithms for control of uncertain systems, see \cite{calafiore_research_2011}, and in particular \cite{alamo_randomized_2012,alamo_randomized_2009,
chamanbaz_sequential_2013, koltchinskii_improved_2000}.
However, we remark that the sample bound~$M_\ic$ in Algorithm \ref{alg:sequential algorithm} is
strictly less conservative than the bound derived in \cite{alamo_randomized_2012} because the infinite sum (Riemann Zeta function) is replaced with a finite sum, following ideas similar to those recently introduced in \cite{chamanbaz_sequential_2013}. This enables us to choose
$\alpha<1$ in (\ref{eq:sample bound Mk}) resulting in up to 30\% improvement in the sample complexity.

Another important difference is on how the cardinality of the design sample set $N_\ic$ appears in the sequential algorithm. In \cite[Algorithm 1]{chamanbaz_sequential_2013}, the constraints are required to be satisfied for all the samples extracted from the set $\mathbb{Q}$ while, in Algorithm \ref{alg:sequential algorithm}, we allow a limited number of  samples to violate the constraints in (\ref{eq:LMI}) and (\ref{eq:BMI}), or their semidefinite versions, in both ``design" and ``validation" steps. Finally, we note that the sequential randomized algorithm in \cite[Algorithm 2]{chamanbaz_sequential_2013} is purely based on additive and multiplicative Chernoff inequalities and hence may provide larger sample complexity than (\ref{eq:sample bound Mk}).

It should also be remarked that the optimal values of the constants $a$ and $\alpha$ depend on other parameters of the algorithm, such as the termination parameter $\ic_t$,  the accuracy level $\varepsilon$, and the level parameter $\rho$. Suboptimal values of $a$ and $\alpha$ for which the sample bound (\ref{eq:sample bound Mk}) is minimized for a wide range of probabilistic levels are $a=3.05$ and $\alpha=0.9$. Note also that for $\rho=0$ the optimal value of these parameters is $a=\infty$ and $\alpha=0.1$, and the bound (\ref{eq:sample bound Mk}) reduces to bound (12) in \cite{chamanbaz_sequential_2013}.
The parameter $\rho$ plays a key role in the algorithm. Note that, as pointed out in \cite{alamo_randomized_2012}, 
the choice of $\rho=0$ may lead to an unfeasible optimization problem in (\ref{eq:optimization in seq alg}) whenever the original robust LMI/BMI is unfeasible.  On the other hand, if $\rho>0$, problem (\ref{eq:optimization in seq alg}) becomes immediately a mixed-integer program, which is numerically hard to solve even in the LMI case.

Finally, we point out that  the termination parameter $\ic_t$ defines the maximum number of iterations of the algorithm which can be chosen by the user. For problems in which the bound $N_{\text{MI}}$
in Algorithm \ref{alg:sequential algorithm}  is large, larger values of $\ic_t$ may be used. In this way, the sequence of sample bounds $N_\ic$ would start from a reasonably small number and would not increase dramatically with the iteration counter $\ic$.
\end{remark}

\section{Numerical Simulations}\label{sec:simulation}
We illustrate the effectiveness of the previous results for a linear model of a robot manipulator joint taken from \cite{kanev_controller_2000}. The state-space model of the plant is given by
\begin{align*}
\begin{cases}
\dot{x}(t)=Ax(t)+Bu(t)\\
z(t)=C_1x(t)+D_{11}w(t)\\
y(t)=Cx(t)+D_{21}w(t)\\
\end{cases}
\end{align*}
where
\begin{align*}
& A=\left[\begin{array}{cccc}
0 & 1 & 0 & 0 \\
0 & 0 &  \frac{c}{M^2I_m} & 0\\
0 & 0 & 0 & 1 \\
0 &  -\frac{\beta}{I_{son}}\,\, &  -\frac{c}{M^2I_m}-\frac{c}{I_{son}}\,\, &  -\frac{\beta}{I_{son}}\\
\end{array}\right],\quad \\
& B=\left[\begin{array}{c}
0\\
 \frac{L_t}{MI_m}\\
0\\
 -\frac{L_t}{MI_m}\\
\end{array}\right],\quad
C=\left[\begin{array}{cccc}
0 & M & 0 & 0\\
1 & 0 & 1 & 0\\
\end{array}\right],\quad D_{21}=\left[\begin{array}{c}
1\\
0\\
\end{array}\right],\\
&C_1=\left[\begin{array}{cccc}
1 & 0 & 1 &0\\
\end{array}\right]\text{ and } D_{11}=1.
\end{align*}
The nominal values of the parameters are as follows: gearbox ratio $M=-260.6$,
motor torque constant $L_t= 0.6$,
damping coefficient $\beta=0.4$,
input axis inertia $I_m=0.001$,
output system inertia $I_{son}=400$,
spring constant $c=130$.
We considered all plant parameters to be uncertain by $15\%$. The objective is to design a static output feedback controller which minimizes the worst case $\mathcal{H}_\infty$ norm of the transfer function from the disturbance channel $w$ to the controlled output $z$. This problem can be formulated in terms of a bilinear matrix inequality \citep{leibfritz_compleib:_2004} of the form
\begin{align}\label{BMI-ex}
&\underset{F,X,\gamma}{\text{minimize }}\gamma\\
& \text{subject to }X\succ 0,\nonumber\\
%%%
& 
\scalebox{0.92}{\mbox{$
\left[\begin{array}{ccc}
A_{F}^TX+XA_{F}\,\, & X(B_1+BFD_{21})&(C_1+D_{12}FC)^T\\
\star & -\gamma &(D_{11}+D_{12}FD_{21})^T\\
\star  &  \star &-\gamma
\end{array}
\right]\prec0$
}}
\nonumber
\end{align}
where $A_{F}\doteq A+BFC$, $X=X^T\in\mathbb{R}^{4\times 4}$, $F\in\mathbb{R}^{1\times2}$ and $\star$ denotes entries that follow from symmetry.

Algorithms \ref{alg:randomized algorithm} and \ref{alg:sequential algorithm} were implemented
using the \textit{Randomized Algorithm Control Toolbox} (RACT) \citep{tremba_ract:_2008}, and we used PENBMI \citep{kocvara_pennon:_2003} for solving BMI optimization problems. The probability density functions of all $6$ uncertain parameters was assumed to be uniform. The level parameter $\rho$ in all simulations was chosen to be zero ($\rho=0$).
A bound on the VC-dimension of the BMI problem (\ref{BMI-ex}) can then be computed using Theorem~\ref{theorem:VC dimension of LMI and BMI}, taking into account that $m_{x}+m_{y}=13$ (for the design variables $F,X$ and $\gamma$), and that $n=6+4+1=11$. Applying Corollary \ref{corollary:one sided sample bound}, the corresponding bounds necessary for applying Algorithm 1 can be computed, and are reported in Table~\ref{tab: simulation results} (third column) for different values of $\delta$ and $\epsilon$.

\begin{table*}[tb]
\begin{center}
\scalebox{0.75}{
%\centering{}%
\begin{tabular}{cc|c||c|ccc|ccc|ccc|ccc}
\toprule
$\varepsilon$ & $\delta$ & \emph{Bound (\ref{eq:sample bound BMI twosided})}  &$\ic_t$ &  \multicolumn{3}{c}{\emph{Design samples}} & \multicolumn{3}{c}{\emph{Validation samples}} & \multicolumn{3}{c}{\emph{Objective value}} & \multicolumn{3}{c}{\emph{Iteration}}
\tabularnewline
%  &  & \emph{Bound}  & \multicolumn{3}{c}{\emph{Samples}} & \multicolumn{3}{c}{\emph{Samples}} &
%\multicolumn{3}{c}{\emph{Value}} & \multicolumn{3}{c}{\emph{Number}}\tabularnewline
\midrule
\midrule
&  &  &  &   \emph{Mean} & \emph{Standard} & \emph{Worst} & \emph{Mean} & \emph{Standard} & \emph{Worst} & \emph{Mean} & \emph{Standard} & \emph{Worst} & \emph{Mean} & \emph{Standard} & \emph{Worst}\tabularnewline
&  &  &  & & \emph{Deviation} & \emph{Case} &  & \emph{Deviation} & \emph{Case} &  & \emph{Deviation} & \emph{Case}&  & \emph{Deviation} & \emph{Case}\tabularnewline
\midrule
$0.2$ & $10^{-2}$ & $3.58\times10^4$ & $5\times10^3$  & $60.6$ & $24.04$ & $149$ & $56.74$ & $0.44$ & $57$ & $1.01$ & $0$ & $1.01$ & $4.8$ & $1.9$& $12$\tabularnewline
\midrule
$0.1$ & $10^{-4}$ & $8.12\times10^4$ & $5\times10^3$ & $149.5$ & $58.7$ & $336$ & $163.2$ & $0.49$ & $164$ & $1.01$ & $0$ & $1.01$ & $5.34$ & $2$ & $12$\tabularnewline
\midrule
$0.05$ & $10^{-6}$ & $1.82\times10^5$ & $10^4$ & $268.7$ & $117.8$ & $594$ & $437.5$ & $0.98$ & $439$ & $1.01$ & $0.01$ & $1.11$ & $8.6$ & $3.7$ & $19$\tabularnewline
\midrule
$0.01$ & $10^{-8}$ & $1.13\times10^6$ & $10^4$  &$1276.5$ & $484.8$ & $2522$ & $2686.5$ & $3.9$ & $2694$ & $1.01$ & $0$ & $1.01$ & $6.6$ & $2.5$ & $13$\tabularnewline
\midrule
$0.005$ & $10^{-10}$ & $2.45\times10^6$ & $10^4$ & $2881.9$ & $1093.3$ & $6310$ & $6305.9$ & $7.9$ & $6323$ & $1.01$ & $0$ & $1.01$ & $6.8$ & $2.6$ & $15$\tabularnewline
\bottomrule
\end{tabular}
}
\end{center}
\vskip .15in
\caption{Sample complexity bounds and simulation results obtained using Algorithm \ref{alg:sequential algorithm}. The third column is the original sample complexity bound (\ref{eq:sample bound BMI twosided}) for strict BMIs, and the fifth column is the sample complexity achieved using Algorithm \ref{alg:sequential algorithm}.}
\label{tab: simulation results}
\end{table*}

Clearly, these sample bounds are quite large. For this reason, we used Algorithm \ref{alg:sequential algorithm} to efficiently solve the problem. Since the sample complexities $M_k$ and $N_k$ in which  Algorithm \ref{alg:sequential algorithm} terminates are random variables, we run the simulations $100$ times for each pair of probabilistic accuracy and confidence parameters. The mean, standard deviation and worst case values of the number of design samples, validation samples, objective value and the iteration number in which the algorithm exits are tabulated in Table \ref{tab: simulation results}.
We conclude that with Algorithm \ref{alg:sequential algorithm} we can achieve the same probabilistic levels with a much smaller number of design samples.

\section{Conclusions}

In this paper, we computed explicit bounds on the Vapnik-Chervonenkis dimension (VC-dimension) of two problems frequently arising in robust control, namely the solution of uncertain LMIs and BMIs. In both cases, we have shown that the VC-dimension is linear in the number of design variables. These bounds are then used in a sequential randomized algorithm that can be efficiently applied to obtain probabilistic optimal solutions to uncertain LMI/BMI. Since the sample complexity is independent of the number of uncertain parameters, the proposed algorithm runs in polynomial time.

\section*{Acknowledgements} The results presented in this paper were obtained when the first author was visiting CNR-IEIIT for a period of six months under DSI funding.

\appendix
\section{Proof of Theorem \ref{theorem:VC dimension of LMI and BMI}}\label{app:proof of vc dim for lmi and bmi}
First, we introduce the following definition.
\begin{definition}[$(\gamma,\eta)$-Boolean Function]\label{def: alpha k}
The function $g:\mathbb{R}^\nt\times\mathbb{Q}\to\{0,1\}$ is a $(\gamma,\eta)$-Boolean function if for fixed $q$ it can be written as expressions consisting of Boolean operators involving $\eta$ polynomials
\[\beta_1(\theta),\ldots,\beta_\eta(\theta)\]
in the components $\theta_i,\,i=1,\ldots,\nt$ and the maximum degree of these polynomials with respect to $\theta_i,\,i=1,\ldots,\nt$ is no larger than $\gamma$.
\end{definition}
The following lemma \citep{vidyasagar_learning_2002} which is an improvement on the original result of \cite{karpinski_polynomial_1995}, states an upper bound on the VC-dimension of $(\gamma,\eta)-$Boolean functions.
\begin{lemma}\label{lemma: vc dim for bolean}
Suppose that $g:\mathbb{R}^\nt\times\mathbb{Q}\to\{0,1\}$ is an $(\gamma,\eta)$-Boolean function, then
\begin{equation}\label{eq:vc dim for boolean}
{\rm VC}_g\leq 2\nt\lg(4e\gamma \eta).
\end{equation}
\end{lemma}
In  view of this lemma, in order to find the VC-dimension of the uncertain LMI and BMI problems, it suffices to represent the constraints in (\ref{eq:LMI}) and (\ref{eq:BMI}) as polynomial inequalities. It is well known that an $n\times n$ real symmetric  matrix is positive definite if and only if all $2^n$ principal minors are positive. However, this condition is equivalent to checking positivity of all $n$ leading principal minors.

Since LMIs are a special case of BMIs, we first prove Theorem \ref{theorem:VC dimension of LMI and BMI} for the more general case of BMIs. Let $F_{{\rm BMI}, ij}(x,y,q)$ be the $ij$-th element of the BMI in (\ref{eq:BMI}). The leading principal minors of $F_{{\rm BMI}}(x,y,q)$ are
\begin{align*}
& F_{{\rm BMI}, 11}(x,y,q), \\
& \det\left(\left[\begin{array}{cc}
F_{{\rm BMI}, 11}(x,y,q) & F_{{\rm BMI}, 12}(x,y,q)\\
F_{{\rm BMI}, 21}(x,y,q)   & F_{{\rm BMI}, 22}(x,y,q)\\
\end{array}\right]\right),\ldots,\\
& \det\left(\left[\begin{array}{ccc}
F_{{\rm BMI}, 11}(x,y,q) & \cdots & F_{{\rm BMI}, 1k}(x,y,q)\\
\vdots & &\vdots\\
F_{{\rm BMI}, k1}(x,y,q)& \cdots & F_{{\rm BMI}, kk}(x,y,q)\\
\end{array}\right]\right),
\ldots, \\
& \det\left(F_{{\rm BMI}}(x,y,q)\right).
\end{align*}
Since the number of leading principal minors is $n$, we need to check $n$ polynomial inequalities.
\noindent Next, we need to find the maximum degree of each polynomial inequality with respect to design variables $x_i,\,i=1,\ldots,\nx$ and $y_j,\,j=i,\ldots,\ny$. Based on the definition of determinant, $k$-th leading principal minor of the BMI in (\ref{eq:BMI}) for $k=3,\ldots,n$ can be written as
\begin{equation}\label{eq:leading principal minor}
D_k=\sum_{\ell=1}^k(-1)^{\ell+1}F_{{\rm BMI}, \ell 1}(x,y,q)M_{\ell 1}
\end{equation}
where $D_k$ is the $k$-th principal minor and $M_{\ell 1}$ is the $(\ell,1)$ minor of a matrix formed by the first $k$ rows and columns of the BMI in (\ref{eq:BMI}).
Then, we have that the $k$-th leading principal minor has  maximum degree $k$ with respect to design variables $x_{i}, i=1,\ldots,\nx$ and $y_{j}, j=1,\ldots,\ny$. From the definition of the BMI in (\ref{eq:BMI}), it is clear that every element of the BMI, including the first leading principal minor, has maximum degree $1$ with respect to the design variables $x_{i}, i=1,\ldots,\nx$ and $y_{j}, j=1,\ldots,\ny$. The second leading principal minor of the BMI in (\ref{eq:BMI})
\begin{align*}
D_2= & F_{{\rm BMI}, 11}(x,y,q) F_{{\rm BMI}, 22}(x,y,q)-\\
&F_{{\rm BMI}, 21}(x,y,q) F_{{\rm BMI}, 12}(x,y,q)
\end{align*}
is a polynomial of maximum degree $2$. For $k>2$, the maximum degree of $D_k$ in (\ref{eq:leading principal minor}) is defined by the multiplication of $F_{{\rm BMI}, \ell 1}(x,y,q)$ and $M_{\ell 1}$. The former has the maximum degree $1$ and the latter has maximum degree equal to $D_{k-1}$ because they are of the same order. Hence, the maximum degree of the  $k$-th leading principal minor with respect to the design variables for $k=1,\ldots,n$ is $k$.

Therefore, checking positive definiteness of the BMI in (\ref{eq:BMI}) is equivalent to checking $n$ polynomial inequalities of degree ranging from $1$ to $n$ which can be represented as an $(\gamma,\eta)-$Boolean function with $\gamma=\eta=n$. The result of Theorem \ref{theorem:VC dimension of LMI and BMI} follows by substituting the obtained values of $\gamma$ and $\eta$ into (\ref{eq:vc dim for boolean}). We notice that the same reasoning holds for the case of LMI and we can represent the LMI in (\ref{eq:LMI}) as an $(\gamma,\eta)-$Boolean function with $\gamma=\eta=n$.

\section{Proof of Theorem \ref{theorem:VC dim semidef}}\label{app:proof of nonstrict}

The result follows observing that an $n\times n$ symmetric matrix is positive semidefinite if and only if \textit{all} $2^n$ principal minors are nonnegative. Then, following similar reasoning as in the proof of Theorem \ref{theorem:VC dimension of LMI and BMI}, it follows that checking positive semidefiniteness of (\ref{eq:LMI}) and (\ref{eq:BMI}) is equivalent to evaluating $2^n$ polynomial inequalities of degree ranging from $1$ to $n$. This can be represented as $(\gamma,\eta)-$Boolean function with $\gamma=n$ and $\eta=2^n$. The results of Theorem \ref{theorem:VC dim semidef} follow by substituting the obtained values of $\gamma$ and $\eta$ in (\ref{eq:vc dim for boolean}).

\bibliographystyle{model5-names}
\small
\bibliography{ref,Frugi-biblio}

\end{document}